\documentclass[a4paper,11pt]{article}   % MikTeX
\usepackage[utf8]{inputenc} % utf8 cp866 cp1251
\usepackage{hyperref}
\usepackage[english]{babel}
\usepackage{ifthen}
\usepackage{amssymb,amsmath,amsthm}

 \newif\ifNoRemark
    \def\addtheorem#1#2#3#4{ % \usepackage{ifthen} needed
    \ifthenelse{\expandafter\isundefined\csname the#2\endcsname}{\newcounter{#2}}{}
    \newenvironment{#1}[1][\global\NoRemarktrue]% No Remark by default
     {\par\addvspace{2mm}\noindent
       \refstepcounter{#2}{\bf #3~\csname the#2\endcsname
      \vphantom{##1}\ifNoRemark.\ \else\ (##1).\fi}\begingroup #4}%
     {\endgroup\par\addvspace{1mm}\global\NoRemarkfalse}
    \expandafter\newcommand\csname b#1\endcsname{\begin{#1}}
    \expandafter\newcommand\csname e#1\endcsname{\end{#1}}
    }

\newtheorem{theorem}{Theorem}

\newtheorem{proposition}{Proposition}
\newtheorem{corollary}{Corollary}

\def\bproof{Proof.\,}
\def\proofend{$\square$\vspace{0.3em}\par}

\begin{document}

\title{On the number of relevant variables for discrete functions}
 \date{}

\author{ V. N. Potapov\\
% independent researcher,
{Novosibirsk, Siberia, vpotapov@math.nsc.ru}}

\maketitle

    \begin{abstract}
We consider various definitions of degrees of discrete functions and
establish relations between the number of relevant (essential)
variables and degrees of two- and three-valued functions.

    \end{abstract}

\section{Introduction}

Let $T$ be an arbitrary set and let $T^n$ be the Cartesian power of
$T$. Given a function $f$ on $T^n$, a variable $x_i$, $1\leq i\leq
n$, is called {\it relevant} (essential, or effective) if there
exist $a_1,\dots,a_{i-1}, a_{i+1},\dots, a_n \in T$ and $b, c\in T$
such that
$$f(a_1,\dots,a_{i-1},b, a_{i+1},\dots, a_n)\neq f(a_1,\dots,a_{i-1},c, a_{i+1},\dots,
a_n).$$

In this paper we study the relationship between various concept of
degrees and the number of relevant variables for two- and
three-valued functions on $[q]^n$, where $[q]$ is a $q$-element set.
Binary-valued functions can be considered  as  indicator functions
of subsets of $[q]^n$, so we can speak about the number of relevant
variables for sets. For any bijections $\pi:T\rightarrow [q]$ and
$\sigma:[p]\rightarrow P$ relevant variables of $f:[q]^n\rightarrow
[p]$ one-to-one correspond to  relevant variables of  $\sigma\circ f
\circ \pi$. So, the relevance of variables does not depend on the
domain and the image sets of functions. For convenience, we take
$\{-1,1\}$ as the image set of two-valued functions. Binary-valued
functions on $\{0,1\}^n$ are called Boolean.

 It is easy to see that every Boolean function $f:\{0,1\}^n
\rightarrow \{-1,1\}$ can be represented as a real polynomial. The
minimum degree of a polynomial that coincides with $f$ on
$\{0,1\}^n$ is called the degree of $f$.

 A famous theorem of Nisan and Szegedy \cite{NisSze} states that
a Boolean function  of degree $d$  has at most $d2^{d-1}$ relevant
variables. A similar bound, in which the degree of a Boolean
function is replaced by the order of correlation immunity, was
proved in \cite{Tar}. This bound was improved to $6.614 \cdot 2^d$
in \cite{CHS}, and then it was further improved to $4.394 \cdot 2^d$
in \cite{Wel}.

It is possible to generalize the definition of the degree to other
discrete functions in different ways, one of them is used by Filmus
and Ihringer \cite{FilIhr}. The precise definition of this degree
will be given in the next section.  The remark at the end of their
paper \cite{FilIhr} and the upper bound for Boolean functions from
\cite{Wel} imply that a two-valued function $f$ on $[q]^n$ of degree
$d$ has at most $4.394 \cdot 2^{\lceil \log_2 q\rceil d}$ relevant
variables. In \cite{Val24} this bound was improved to
$\frac{dq^{d+1}}{4(q-1)}$ for $q\neq 2^s$.

In the next section we introduce  degrees $\deg_i(f)$  for functions
$f:[q]^n\rightarrow [p]$, where $\deg_0(f)$ coincides with the
degree $d$. We prove  upper bounds $\frac14\pi^2{\rm
deg}_1(f)q^{\deg_0(f)-1}$ and $\frac12\pi^2\deg_2(f)q^{{\rm
deg}_0(f)-2}$ (Theorem 1) for the number of relevant variables of
two-valued functions. Unlike the previous bounds, the new bounds
depend on $\deg_1(f)$ and $\deg_2(f)$. It occurs that they are
better than the previous bounds for some classes of functions. For
example, the second bound is better than others if $q\geq 4$ and
${\rm deg}_2(f)=\deg_0(f)$. Moreover, we obtain upper bounds
$\frac{{\rm deg}_0(f)q^{\deg_0(f)+1}}{3(q-1)}$ (Theorem 2),
$\frac{\pi^2}{3}{\rm deg}_1(f)q^{\deg_0(f)-1}$, and
$\frac{2\pi^2}{3}{\rm deg}_2(f)q^{\deg_0(f)-2}$ (Theorem 3) for the
number of relevant variables in the case of three-valued functions.

Our proofs are based on the notion of an average sensitivity. We
consider a function $f:[q]^n\rightarrow [p]$ as a $p$-coloring of a
graph $G$ such that $|V(G)|=q^n$. The average sensitivity  $I[f]$ is
the number of mixed colored edges in $G$. Our estimation of $I[f]$
is similar to the proof of the Bierbrauer--Friedman bound (see
\cite{Bier} and \cite{Pot}) and depends on the adjacency  matrix of
$G$. In  previous papers \cite{NisSze}, \cite{FilIhr}, \cite{Wel},
\cite{Val24} the authors implicitly or explicitly treated $G$ as the
Hamming graph. In the present paper we use the Cartesian products of
cycles instead of the Hamming graphs.

Moreover, in Section 3 we discuss relations between these degrees
and other well-known degrees of Boolean function such as numerical
and algebraic degrees.

\section{Fourier--Hadamard transform}

In this section we treat the domain $[q]^n$ of functions as an
abelian group $G$ of order $[q]^n$. Consider  the vector space
$V(G)$ consisting of functions $f:G\rightarrow \mathbb{C}$ with the
inner product
$$(f,g)=\sum\limits_{x\in G}f(x)\overline{g(x)}.$$
A function $f:G\rightarrow \mathbb{C}\backslash\{0\}$ mapping from
$G$ to the non-zero complex numbers is called a character of $G$ if
it is a group homomorphism from $G$ to $\mathbb{C}$, i.e.,
$\phi(x+y)=\phi(x)\phi(y)$ for each $x,y\in G$. The set of {\it
characters} of $G$ is an orthogonal basis of $V(G)$.

We  consider the linear space $V(\mathbb{Z}^n_q)$ of complex valued
functions with finite domain
$\mathbb{Z}^n_q=(\mathbb{Z}/q\mathbb{Z})^n$. Let $\xi=e^{2\pi i/q}$.
We can define characters of $\mathbb{Z}^n_q$ as
$\phi_z(x)=\xi^{\langle x,z\rangle}$, where  and $\langle
x,z\rangle=x_1z_1+\dots+x_nz_n\, {\rm mod}\, q$ for each $z\in
\mathbb{Z}^n_q$.

Below we will consider $\mathbb{Z}_q$ as the set
$\{-\frac{q-2}{2},\dots,-1,0,1,\dots,q/2\}$ if $q$ is even and  as
the set $\{-\frac{q-1}{2},\dots,-1,0,1,\dots,\frac{q-1}{2}\}$ if $q$
is odd. We define the $m$th degree of $\phi_z$, $z=(z_1,\dots,z_n)$,
as the sum $\deg_m(\phi_z)=\sum_{k=1}^{n}|z_k|^m$. A {\it weight} of
$z\in \mathbb{Z}^n_q$ is the number of nonzero coordinates of $z$,
i.\,e., ${\rm wt}(z)=\deg_0(\phi_z)$.

Changing the variables $x_i\rightarrow y_i=\xi^{x_i}$ or
$x_i\rightarrow y_i=\xi^{-x_i}$ we see that  $\phi_z$ corresponds to
an ordinary monomial of degree $\deg_1\phi_z$.

Consider the expansion of $f\in V(\mathbb{Z}^n_q)$ with respect to
the basis of characters
\begin{equation}\label{vareq0}
f(x)=\frac{1}{q^{n}}\sum\limits_{z\in \mathbb{Z}^n_q}
W_f(z)\phi_z(x), \end{equation} where $W_f(z)=(f,\phi_z)$ are called
the {\it Fourier--Hadamard coefficients} of $f$. The function
$W_f\in V(\mathbb{Z}^n_q)$  is called the {\it Fourier--Hadamard} or
{\it Walsh--Hadamard} (in binary case) {\it transform}  of $f$. We
define
$$\deg_m(f)=\max\limits_{W_f(z)\neq 0} \deg_m(\phi_z).$$

If $q=2$ or $q=3$ then we see that  $\deg_m(f)=\deg_0(f)$ for all
$m$. Note that in \cite{FilIhr} and \cite{Val24} the authors call
$\deg_0(f)$  a degree of $f$.

\section{Properties of numerical degree of Boolean functions}

 Let
$T$ be a finite subset of $\mathbb{C}$. Consider the linear space
$V(T^n)$ of complex valued functions on $T^n$. Let
$C_k(x_1,\dots,x_n)$ be the linear space of polynomials over
$\mathbb{C}$, where every variable has degree at most  $k-1$.
\begin{proposition}\label{valprop100}
For every  function $f\in V(T^n)$  there exists unique polynomial
$P_f\in C_k(x_1,\dots,x_n)$, $k=|T|$, such that $P_f|_{T^n}=f$.
\end{proposition}
Proof. We will prove the existence of the polynomial by induction.
 If $n=1$ then $P_f$ is the Lagrange interpolating polynomial. By the induction
hypothesis, there exist $P_i|_{T^{n-1}\times
\{t_i\}}=f|_{T^{n-1}\times \{t_i\}}$, where $t_i\in T$. Then
$P_f(\overline{x})=\sum\limits_{i=1}^nP_i(\widetilde{x}_i)\frac{\prod\limits_{t_j\in
T\setminus\{t_i\}}(x_i-t_j)}{\prod\limits_{t_j\in
T\setminus\{t_i\}}(t_i-t_j)}$, where $\widetilde{x}_i$ is the set of
all variables except for $x_i$. Since the dimensions of $V(T^n)$ and
$C_k(x_1,\dots,x_n)$ coincide, such a polynomial is unique.
\proofend

Let $\deg P$ be the degree of $P\in C_k(x_1,\dots,x_n)$ and let
$\deg' P$ be the maximum number of variables in the monomials of
$P$. Obviously, $\deg' P\leq \deg P$ and if $k=2$ then $\deg' P=\deg
P$. We define $\deg_{num}f=\deg P_f$ and $\deg'_{num}f=\deg' P_f$.

\begin{proposition}\label{valprop101}
Let $s:\mathbb{Z}_q\rightarrow \mathbb{C}$ be defined by the
equation $s(x)=\xi^x$. Suppose that $f=g\circ s$, where
$f:\mathbb{Z}_q^n\rightarrow \mathbb{C}$ and $g\in
V((s(\mathbb{Z}_q))^n)$. Then $\deg'_{num}g=\deg_0f$ and
$\deg_{num}g\geq \deg_1f$.
\end{proposition}
Proof. By (\ref{vareq0}) we have
$$g(x_1,\dots,x_n)=\frac{1}{q^{n}}\sum\limits_{z\in \mathbb{Z}^n_q}
W_f(z)x_1^{z_1}\cdots x_n^{z_n}=\frac{1}{q^{n}}\sum\limits_{y\in
\{0,1,\dots,q-1\}^n} W_f(z)x_1^{y_1}\cdots x_n^{y_n},$$ where
$y_i=z_i\mod q$. Therefore, $\deg_{num}g=\max_y\sum_{k=1}^{n}y_k\geq
\max_z\sum_{k=1}^{n}|z_k|=\deg_1f$. Moreover,
$\deg'_{num}g=\max\limits_{W_f(z)\neq 0} {\rm wt}(z)=\deg_0f$.
\proofend

Consider a function  $f:\mathbb{Z}^n_q\rightarrow \mathbb{C}$ and
two surjections $s_i:\mathbb{Z}_q\rightarrow \mathbb{C}$, $i=1,2$.
Let $T_i=s_i(\mathbb{Z}_q)$ and $f=g_i\circ s_i$, $i=1,2$. Then
$g_i\in V(T_i^n)$, $i=1,2$. It is easy to see that $\deg'_{num}g_1=
\deg'_{num}g_2$ but $\deg_{num}g_1$ and $\deg_{num}g_2$ may be
different even in the case $n=1$. So if we want to define the degree
$\deg_{num}f$ as $\deg_{num}g_1$ then it  will unfortunately depend
on the surjection of a finite set into $\mathbb{C}$. Below we will
consider the case $|T|=2$ in  more detail. In this case
$\deg_{num}f=\deg'_{num}f$ and therefore  this degree does not
depend on the surjection into $\mathbb{C}$.

Next we treat the domain $[2]^n$ of functions as a vector space
$\mathbb{F}_2^n$. A real valued function  $f: \mathbb{F}_2^n
\rightarrow \mathbb{R}$ is called a {\it pseudo-Boolean function}.
By Proposition \ref{valprop100} every pseudo-Boolean function can be
represented in  {\it numerical normal form} (NNF)
\begin{equation}\label{eqZhegal1}
f(x_1,\dots,x_n)=\sum\limits_{y\in
    \mathbb{F}_2^n}a(y)x_1^{y_1}\cdots x_n^{y_n},
\end{equation}
     where $x^0=1, x^1=x$, and $a(y)\in\mathbb{R}$.
The  maximum degree of the monomial in NNF is called the {\it
numerical degree} of $f$.
%, i.\,e., $
%\deg_{num}(f)=\max\limits_{a(y)\neq 0} {\rm wt}(y)$ by Proposition
%\ref{valprop101}.

 Every Boolean function $f: \mathbb{F}_2^n
\rightarrow \mathbb{F}_2$ can be represented in  {\it algebraic
normal form} (ANF)
\begin{equation}\label{eqZhegal}
f(x_1,\dots,x_n)=\bigoplus\limits_{y\in
    \mathbb{F}_2^n}M_f(y)x_1^{y_1}\cdots x_n^{y_n},
\end{equation}
     where $x^0=1, x^1=x$,
and the function $M_f:\mathbb{F}_2^n\rightarrow \mathbb{F}_2$ is
called the {\it M\"obius transform} of $f$. It is well known (see
\cite{Carlet}) that for every function ANF is unique.
 %(see, e.g. \cite{Carlet} Section 2.2.1).

 The  maximal degree of the monomial in ANF of
$f$ is called {\it algebraic degree} of $f$, i.\,e., ${
\deg}_{alg}(f)=\max\limits_{M_f(y)=1} {\rm wt}(y)$. If all monomials
have degree one, then $f$ is a {\it linear} function. Denote by
$\ell_u$ the linear function $\ell_u(x)=\langle u,x
\rangle=u_1x_1\oplus u_2x_1\oplus\dots \oplus u_nx_n$, where $u\in
\mathbb{F}^n_2$, and $\ell_{\bf 1}(x)=x_1\oplus x_2\oplus\cdots
\oplus x_n$. Obviously, if $f\neq const$ or $f\neq \ell_{\bf 1}$
then $\deg_{alg}(f)=\deg_{alg}(f\oplus \ell_{\bf 1})$. A variable
$x_i$ of $f$ is called {\it linear} if $f(x_1,\dots,x_n)=
g(x_1,\dots,x_{i-1}, x_{i+1},\dots,x_n)\oplus x_i$.

We can consider a Boolean function as a pseudo-Boolean function with
values $\{0,1\}\subset \mathbb{R}$. It is easy to prove that
$\deg_{alg}(f)\leq \deg_{num}(f)$ for  any Boolean function $f$.
Indeed, consider ANF $f(x_1,\dots,x_n)=\bigoplus\limits_{y\in
    \mathbb{F}_2^n}a(y)x_1^{y_1}\cdots x_n^{y_n}$, where
    $a(y)=M_f(y)$. Then

$(-1)^{f(x_1,\dots,x_n)}=\prod\limits_{y\in
    \mathbb{F}_2^n}(-1)^{a(y)x_1^{y_1}\cdots x_n^{y_n}}$, and
$1-2f(x)=\prod\limits_{y\in
    \mathbb{F}_2^n}(1-2a(y)x_1^{y_1}\cdots x_n^{y_n})$,
    since $(-1)^b=1-2b$ for $b\in \{0,1\}\subset \mathbb{R}$.

Using equality $x^2=x$ for $x\in \{0,1\}\subset \mathbb{R}$, we
obtain that $$\deg_{alg}(f)\leq \deg_{num}(f)=\deg_{num}((-1)^f).$$

Denote by $V(\mathbb{F}^n_2)$ the $2^n$-dimensional vector space
(over $\mathbb{R}$) of pseudo-Boolean functions.  By (\ref{vareq0}),
we have

$$(-1)^f(x) = \frac{1}{2^{n}}\sum\limits_{y \in \mathbb{F}_2^n} W_f(y)(-1)^{ \langle y,x \rangle},$$
where  $W_f(y)$ are  the  Walsh--Hadamard  coefficients of $f$.
Since $(-1)^{ \langle y,x
\rangle}=\prod_{i=1}^n(-1)^{y_ix_i}=\prod_{i=1}^n(1-2y_ix_i)$, we
have
$$(-1)^f(x) = \frac{1}{2^n}\sum\limits_{y \in \mathbb{F}_2^n} W_f(y)\prod_{i=1}^n(1-2y_ix_i).$$
Then
\begin{equation}\label{e1variab0l}
\deg_{num}(f)=\deg_{num}((-1)^f)=\max\limits_{W_f(y)\neq 0} {\rm
wt}(y)=\deg_0(f).
\end{equation}

\begin{proposition}
For every Boolean function $f$  it holds $\deg_{alg}(f)\leq
\min\{\deg_0(f),n-\deg_0(f)\}$.
\end{proposition}
Proof. Denote by $\mathcal{W}(f)$ the multiset of  Walsh--Hadamard
coefficients of $f$.  From the definitions we see that
\begin{equation}\label{e1variabl}
y\in \mathcal{W}(f) \Leftrightarrow y\oplus{\bf 1} \in
\mathcal{W}(f\oplus \ell_{\bf 1}). \end{equation}
  Then
$\deg_{num}(f\oplus \ell_{\bf 1} )= n- \min\limits_{W_f(y)\neq 0}
{\rm wt}(y)$. Since $\deg_{alg}(f)=\deg_{alg}(f\oplus \ell_{\bf 1})$
if $\deg_{alg}(f)>1$, then we obtain another inequality
$\deg_{alg}(f)\leq \min\{\max\limits_{W_f(y)\neq 0} {\rm wt}(y), n-
\min\limits_{W_f(y)\neq 0} {\rm wt}(y)\}$. By (\ref{e1variab0l}) we
obtain the required inequality if $\deg_{alg}(f)>1$. For
$\deg_{alg}(f)\leq1$ the required inequality  is obviously true.
\proofend

Denote by $t(f)$  the number of the relevant variables of $f$. From
the definitions, we have $\deg_{alg}(f)\leq t(f)$ for  Boolean and
$\deg_{num}(f)\leq t(f)$ for pseudo-Boolean functions. Does there
exist a reversed inequality in a general case? There exists a
Boolean function $\ell_{\bf 1}$ with minimal algebraic degree
$\deg_{alg}(\ell_{\bf 1})=1$ and maximal number $n$ of the relevant
variables. Moreover, there exists a pseudo-Boolean function
$\jmath(x)=(-1)^{x_1}+\cdots +(-1)^{x_n}= n-2(x_1+\cdots+x_n)$ with
minimal numerical degree $\deg_{num}(\jmath)=1$  and maximal number
$n$ of the relevant variables. Thus, the  inequalities for algebraic
degree of Boolean functions and for numerical degree of
pseudo-Boolean functions cannot be  reversed. However, as mentioned
in  Introduction, the numerical degree provides an upper bound for
the number  of relevant variables in the case of Boolean functions.
In the next sections we prove upper bounds for the number of the
relevant variables for $q$-ary two- and three-valued functions.

\section{Bounds for two-valued functions}

The {\it Cayley graph} $Cay(G, S)$ on abelian group $G$ with
connecting set $S$, $S\subset G$, $S = -S$, $0\not \in S$, is the
graph whose vertices are the elements of $G$ and whose edge set $E$
is $\{\{x, a+x\} : x\in G, a\in S\}$.

It is well known that the set of scalar characters of abelian group
$G$ is an orthogonal basis consisting of the eigenfunctions of
$Cay(G, S)$. The eigenfunctions of a graph $\Gamma$ are eigenvectors
of the adjacency matrix of $\Gamma$.

\begin{proposition}[\cite{Bab}, Corollary 3.2]\label{varprop1}
Let $\phi$ be a character of $\mathbb{Z}^n_q$. Then its eigenvalue
with respect to $Cay(\mathbb{Z}^n_q, S)$ is equal to
$\sum\limits_{s\in S}\phi(s)$.
\end{proposition}

Let $S\subseteq \mathbb{Z}_q\setminus \{0\}$. Consider
$S^n=\{(0,\dots,0,\underset{i}{s},0,\dots,0): s\in S,
i=1,\dots,n\}\subset \mathbb{Z}_q^n$ as a connecting set in
$\mathbb{Z}_q^n$. If $S=\mathbb{Z}_q\setminus \{0\}$ then
$Cay(\mathbb{Z}_q, S)$ is the complete graph $K_q$.  By the
definition of the Cayley graph we obtain that $Cay(\mathbb{Z}_q^n,
S^n)=K_q\square\cdots\square K_q$. This graph is equal to the
Hamming graph $H(n,q)$. The Hamming graph induces the Hamming
distance $d_H$ between vertices. This distance $d_H(u,v)$ is equal
to the number of places in which $n$-tuples $u,v\in \mathbb{Z}_q^n$
differ. The eigenvalues of the Hamming graphs are well known and are
obtained from Proposition \ref{varprop1}.

\begin{corollary}\label{varcor01}
The eigenfunction $\phi_z(x)=\xi^{\langle x,z\rangle}$ corresponds
to the eigenvalue \\ $\lambda_z=(q-1)n- q{\rm wt}(z)$ in $H(n,q)$.
\end{corollary}

In the present paper we take $S=\{-1,1\}$. Thus, $Cay(\mathbb{Z}_q,
S)$ is the circular graph $C_q$ consisting of one cycle. In this
case $S^n$ is a collection of $n$-dimensional vectors consisting of
$\pm 1$ and zeros. Then $Cay(\mathbb{Z}_q^n,
S^n)=C_q\square\cdots\square C_q=C^n_q$. The graph $C^n_q$ is called
a {\it hypercube with induced Lee distance} $d_L$, where
$d_L(u,v)=\sum_{i=1}^n\min\{|u_i-v_i|, q-|u_i-v_i|\}$. If $q=2$ or
$q=3$, then the Hamming and Lee distances are the same. We  say that
an edge $\{x,y\}$ in $C^n_q$  has direction $i$ if vertices $x$ and
$y$ differ in the $i$th position.

For given a vector $z\in\mathbb{Z}^n_q$ denote by $a_k(z)$ the
number of elements $k\in \mathbb{Z}_q$ in $z$. Then the { weight} of
$z$ is ${\rm wt}(z)=\sum\limits_{k\in \mathbb{Z}_q}a_k(z)$.  By
Proposition \ref{varprop1}, we obtain

\begin{corollary}\label{varcor1}
The eigenfunction $\phi_z(x)=\xi^{\langle x,z\rangle}$ corresponds
to the eigenvalue \\ $\lambda_z=2n- 4\sum\limits_{k\in
\mathbb{Z}_q}a_k(z)\sin^2\frac{\pi k}{q} $ in $C^n_q$.
\end{corollary}
\bproof
$$\lambda_z=\sum\limits_{i=1}^n(\xi^{z_i}+\xi^{-z_i}) =
\sum\limits_{i=1}^n 2\cos\frac{2\pi z_i}{q}$$ $$=2n-
\sum\limits_{k\in \mathbb{Z}_q}2a_k(z)(1-\cos\frac{2\pi k}{q})=
2n-\sum\limits_{k\in \mathbb{Z}_q}4a_k(z)\sin^2\frac{\pi k}{q}.$$
\proofend

We will use some results on the theory of invariant subspaces of
Hamming graphs developed by Valyuzhenich and his coauthors. Denote
by $U_{k}(n,q)$ the linear span of all $\phi_z$, where $z$ has
weight $k$. $U_{k}(n,q)$ is a subspace of $V(\mathbb{Z}_q)$. The
direct sum of subspaces
$$U_{0}(n,q)\oplus\cdots \oplus U_{m}(n,q)$$ is denoted by
 $U_{[0,m]}(n,q)$. Straightforwardly, $U_{[0,m]}(n,q)$ is the set of
 functions $f$ such that $\deg_0(f)\leq m$ (see \cite{Val24}).

\begin{proposition}[\cite{Val24}, Theorem 1]\label{propVal1}
Let $f\in U_{[0,m]}(n,q)$, where $q\geq 3$ and $f\neq 0$. Then
$|{\rm supp}(f)|\geq q^{n-m}$.
\end{proposition}

Denote by $f|_{x_i=a}$ a retract of $f$, i.\,e.,
$$f|_{x_i=a}(x_1,\dots,x_{i-1}, x_{i+1},\dots, x_n)=
f(x_1,\dots,x_{i-1},a, x_{i+1},\dots, x_n).$$

\begin{proposition}[\cite{Val19}, Lemma 4]\label{propVal11}
If $f\in U_{[0,m]}(n,q)$, $m>0$. Then the difference\\
$f|_{x_i=a}-f|_{x_i=b} $ belongs to $ U_{[0,m-1]}(n-1,q)$.
\end{proposition}

\begin{corollary}\label{propVal111}
If $f|_{x_i=a}\neq f|_{x_i=b} $ then $|{\rm
supp}(f|_{x_i=a}-f|_{x_i=b})|\geq q^{n-\deg_0(f)}$.
\end{corollary}

The next property follows from  the definition of the
Fourier--Hadamard coefficients.

\begin{proposition}\label{propVal00}
If a function  $f$ does not essentially depend on variable $x_i$ and
$z_i\neq 0$ then $W_f(z)=0$.
\end{proposition}

By the definition of degree we obtain
\begin{corollary}
If a function  $f$ has $m$ relevant variables then ${\rm
deg}_0(f)\leq n-m$.
\end{corollary}

Next, we prove the converse statement on the bound of the number of
relevant  variables under conditions on the degrees of functions.
The proof of the following theorem is similar to the arguments from
\cite{Val24} but we use the hypercube with the Lee metric instead of
one with  the Hamming metric.

\begin{theorem}\label{varth}
For a Boolean valued function $f$ on $\mathbb{Z}^n_q$ it holds
$$t(f)\leq \frac{\pi^2}{4}\deg_1(f)q^{\deg_0(f)-1}\qquad \mbox{\rm and}\qquad t(f)\leq \frac{\pi^2}{2}\deg_2(f)q^{{\rm
deg}_0(f)-2},$$ where $t(f)$ is the number of relevant variables of
$f$.
\end{theorem}
\bproof We will consider the domain of $f$ as the vertex set of
$C_q^n$. Let $A$ be the adjacency matrix of $C_q^n$. An edge
$\{x,y\}$ of $C_q^n$ is called mixed colored if $f(x)\neq f(y)$. The
total number of edges of $C_q^n$ is $nq^n$. Denote by $I[f]$ the
number of mixed colored edges of $C_q^n$. Note that the average
number $\frac{I[f]}{|V(H(n,q)|}$ of mixed colored edges   in the
Hamming graph is called the average sensitivity of $f$. But $I[f]$
may be less than the sensitivity of $f$ in the case of $C_q^n$.
Straightforwardly, we can prove that
$$ -(Af,f)=2I[f]-(2nq^n-2I[f]).$$ By the definition of characters, we
obtain that $f=\frac{1}{q^{n}}\sum\limits_{z\in \mathbb{Z}^n_q}
W_f(z)\phi_z,$ and
\begin{equation}\label{valeq01}
(Af,f)= \frac{1}{q^{2n}}\sum\limits_{z\in \mathbb{Z}^n_q}
\lambda_z|W_f(z)|^2(\phi_z,\phi_z).
\end{equation}
It is clear that $(\phi_z,\phi_z)=q^n$. By Corollary \ref{varcor1},
we obtain that

\begin{equation}\label{valeq1}
I[f]=\frac{1}{q^{n}}\sum\limits_{z\in
\mathbb{Z}^n_q}|W_f(z)|^2\sum\limits_{k\in
\mathbb{Z}_q}a_k(z)\sin^2\frac{\pi k}{q}.
\end{equation}
Using $\sin^2y=\sin^2(\pi-y)$ and $\sin^2y\leq y^2$, we have

\begin{equation}\label{valeq2}
I[f]\leq\frac{1}{q^{n}}\sum\limits_{z\in
\mathbb{Z}^n_q}|W_f(z)|^2\sum\limits_{k=-k'_1}^{k_1}a_k(z)\left(\frac{\pi
k}{q}\right)^2,
\end{equation}
where $k_1=\frac{q}{2}$, $k'_1=\frac{q}{2}-1$ if $q$ is even and
$k'_1=k_1=(q-1)/2$ if $q$ is odd. By the definition of degrees,  we
obtain that $\sum\limits_{k=-k'_1}^{k_1}a_k(z)k^2\leq {\rm
deg}_2(f)$ and $\sum\limits_{k=-k'_1}^{k_1}a_k(z)k^2\leq k_1{\rm
deg}_1(f)$ for all $z\in \mathbb{Z}^n_q$. Then from Parseval's
identity $\sum\limits_{z\in \mathbb{Z}^n_q}|W_f(z)|^2=q^{2n}$ and
(\ref{valeq2}) we obtain
\begin{equation}\label{valeq3}
I[f]\leq  \frac{\deg_2(f)}{q^n}\left(\frac{\pi
}{q}\right)^2\sum\limits_{z\in \mathbb{Z}^n_q}|W_f(z)|^2\leq  {\rm
deg}_2(f)\pi^2q^{n-2}\quad \mbox{\rm and}\quad I[f]\leq 2 {\rm
deg}_1(f)\pi^2 q^{n-1}.
\end{equation}
Let $x_i$ be a relevant variable of $f$. Consider the retracts
$f|_{x_i=0}$, $f|_{x_i=1}$,... There are  at least two numbers $a_1,
a_2\in \mathbb{Z}_q$ such that $f|_{x_i=a_j}\neq f|_{x_i=a_j+1 \mod
q}$, $j=1,2$. By Corollary \ref{propVal111}, we obtain that at least
$2q^{n-\deg_0(f)}$ mixed colored edges have direction $i$. Then
$I[f]\geq 2t(f)q^{n-\deg_0(f)}$. By inequalities (\ref{valeq3}) the
proof is complete.
 \proofend

Next we consider an example of a function $f_m$ such that the new
estimate of $t(f_m)$ is greater than the previous one.
 For $q=3$ the presented bound  $\frac{\pi^2}{2}{\rm
deg}_2(f)q^{d-2}$ is weaker than Valyuzhenich's bound
$\frac{dq^{d+1}}{4(q-1)}$ since $\deg_2(f)\geq \deg_0 (f)=d$ and
$\frac{\pi^2}{2}\geq \frac{3^3}{8}$. So, consider the following
example for $q=4$.
 Let $h:\mathbb{Z}_4\rightarrow \{0,1\}$ be
defined by the vector of values $(1,1,0,0)$.  We have equalities
$\sum_{x\in \mathbb{Z}_4} h(x)i^{-2x}=\sum_{x\in \mathbb{Z}_4}
h(x)i^{2x}=0$, where $i=\sqrt{-1}$. Consider
$f_m:\mathbb{Z}^n_4\rightarrow \{0,1\}$, where
$f_m(x_1,\dots,x_n)=h(x_1)\cdot h(x_2)\cdots h(x_m)$. It is clear
that $t(f_m)=m$. Let us  estimate $t(f_m)$ using the above formulas.
 By Proposition \ref{propVal00}, we conclude that  $W_{f_m}(z)=0$ if
$z_k\neq 0$ for some $k>m$. If $z_k=0$ for all $k>m$, then we obtain
that
$$W_{f_m}(z)=\sum\limits_{x}f_m(x)\xi^{-\langle x,z\rangle} =\sum\limits_{x}f_m(x) \xi^{-\langle x,z\rangle} =\sum\limits_{x}h(x_1)i^{-x_1z_1}\cdots h(x_m)i^{-x_mz_m}$$ $$=4^{n-m}
(\sum_{x_1}h(x_1)\xi^{-x_1z_1})\cdots
(\sum_{x_m}h(x_m)\xi^{-x_mz_m}), \ {\rm where}\ \xi=i.$$ Since
$\sum_x h(x)i^{-xz}=0$ for $z=2$, we conclude that   $ {\rm deg}_2
(f_m)=\deg_0 (f_m)=m$. Thus, the new bound  $t(f_m)\leq
\frac{\pi^2m}{32}4^{m}$ is slightly better than Valyuzhenich's bound
$t(f_m)\leq \frac{m4^{m}}{3}$.

\section{Bounds for three-value functions}

It is possible to generalize our methods to functions with three
different values. We put the set of values $\Xi=\{1,\xi,\xi^{-1}\}$,
where $\xi=e^{\frac{2\pi i}{3}}$.  Let the domain of $f$ be the
vertex set of $C_q^n$ and let $A$ be the adjacency matrix of
$C_q^n$. It is easy to see that $a\overline{b}+\overline{a}b=-1$ if
$a,b\in \Xi$ and $a\neq b$; $a\overline{a}+\overline{a}a=2$ for each
$a\in \Xi$. Then
\begin{equation}\label{valeq151}
(Af,f)= -I[f]+ 2(nq^n-I[f]),
\end{equation}
 where $I[f]$ is the number of mixed
colored edges. Indeed, on the left side of the equation  two
adjacent vertices with equal values give the term $2$ and two
adjacent vertices with different values give the term $-1$.

By (\ref{valeq01}), (\ref{valeq151}) and  Corollary \ref{varcor1} we
obtain that

\begin{equation}\label{valeq15}
I[f]=\frac{4}{3q^{n}}\sum\limits_{z\in
\mathbb{Z}^n_q}|W_f(z)|^2\sum\limits_{k\in
\mathbb{Z}_q}a_k(z)\sin^2\frac{\pi k}{q}.
\end{equation}

Using (\ref{valeq15}) instead of (\ref{valeq1}), similarly to
Theorem \ref{varth} we prove the following inequalities for
three-valued functions.

\begin{theorem}
For a three-valued function $f$ on $\mathbb{Z}^n_q$ it holds
$$t(f)\leq \frac{\pi^2}{3}\deg_1(f)q^{\deg_0(f)-1}\qquad \mbox{\rm and}\qquad t(f)\leq \frac{2\pi^2}{3}\deg_2(f)q^{{\rm
deg}_0(f)-2},$$ where $t(f)$ is the number of relevant variables of
$f$.
\end{theorem}

Moreover, using arguments from \cite{Val24} we can prove the
following statement.

\begin{theorem}
Every three-valued function $f$ of degree $d=\deg_0(f)$ on
$\mathbb{Z}^n_q$,  has at most $\frac{dq^{d+1}}{3(q-1)}$ relevant
variables.
\end{theorem}
\bproof Every vertex of the Hamming graph $H(n,q)$  has $n(q-1)$
neighbors instead of $2n$ neighbors in $C^n_q$. So, if $A$ is the
adjacency matrix of $H(n,q)$, then

\begin{equation}\label{valeq152}
(Af,f)= -I[f]+ 2\left(\frac{n(q-1)}{2}q^n-I[f]\right),
\end{equation}

 By (\ref{valeq01}), (\ref{valeq152}) and  Corollary
\ref{varcor01} we obtain that

\begin{equation}\label{valeq115}
3I[f]=\frac{q}{q^n}\sum\limits_{z\in \mathbb{Z}^n_q}|W_f(z)|^2{\rm
wt}(z).
\end{equation}

Using Parseval's identity $\sum\limits_{z\in
\mathbb{Z}^n_q}|W_f(z)|^2=q^{2n}$ and the definition
$\max\limits_{z}{\rm wt}(z)=\deg_0(f)=d$ we obtain that
\begin{equation}\label{valeq33}
3I[f]\leq q^{n+1}d.
\end{equation}

Let $x_i$ be a relevant variable of $f$. By the definition of the
relevant variable, not all retracts $f|_{x_i=0}$, $f|_{x_i=1}$,...
are equal. Let us estimate the number of pairs of distinct retracts.
Suppose that $t_j$ be the number of  retracts of  type $j$, where
$j=1,\dots,k$, $2\leq k\leq q$, $\sum_{j=1}^k t_j=q$. It is easy to
see that there exist $\sum_{j=1}^k t_j(q-t_j)\geq 2q-2$  ordered
pairs of distinct retracts. Thus, by Corollary \ref{propVal111}, we
obtain that at least $(q-1)q^{n-d}$ mixed colored edges have
direction $i$. Then $I[f]\geq (q-1)t(f)q^{n-d}$. By inequalities
(\ref{valeq33}), the proof is complete.
 \proofend

\end{document}